\documentclass[letterpaper, 10 pt, conference]{ieeeconf}

\IEEEoverridecommandlockouts

\usepackage{amsmath}
\usepackage{amssymb}
\usepackage{xcolor}

\usepackage{subcaption}
\usepackage{cite}
\usepackage{hyperref}
\hypersetup{colorlinks=true,
	urlcolor=blue,
	linkcolor=blue,
	citecolor=blue}

\usepackage{amsthm}

\usepackage{nccmath,mathtools,mathrsfs}

\let\labelindent\relax
\usepackage{enumitem}

\allowdisplaybreaks

\newtheorem{theorem}{Theorem}
\newtheorem{proposition}{Proposition}
\newtheorem{lemma}{Lemma}
\newtheorem{corollary}{Corollary}

\newtheorem{assumption}{Assumption}

\theoremstyle{definition}

\theoremstyle{remark}
\newtheorem{remark}{Remark}

\newcommand{\tr}{{{\mathsf T}}}
\newcommand{\mK}{{\mathsf{K}}}

\usepackage[nameinlink]{cleveref}
\crefname{equation}{}{}
\crefname{theorem}{Theorem}{Theorems}
\crefname{corollary}{Corollary}{Corollaries}
\crefname{example}{Example}{Examples}
\crefname{assumption}{Assumption}{Assumptions}
\crefname{lemma}{Lemma}{Lemmas}
\crefname{proposition}{Proposition}{Propositions}
\crefname{figure}{Figure}{Figures}
\crefname{table}{Table}{Tables}
\crefname{section}{Section}{Sections}
\crefname{appendix}{Appendix}{Appendices}
\Crefname{equation}{}{}
\Crefname{theorem}{Theorem}{Theorems}
\Crefname{corollary}{Corollary}{Corollaries}
\Crefname{example}{Example}{Examples}
\Crefname{lemma}{Lemma}{Lemma}
\Crefname{proposition}{Proposition}{Proposition}
\Crefname{figure}{Figure}{Figures}
\Crefname{table}{Table}{Tables}
\Crefname{section}{Section}{Sections}
\Crefname{appendix}{Appendix}{Appendices}

\title{\LARGE \bf
On the Global Optimality of Direct Policy Search for Nonsmooth $\mathcal{H}_\infty$ Output-Feedback Control
}

\author{Yujie Tang$^{1}$ and Yang Zheng$^{2}$%
\thanks{The work of Y. Zheng is supported by NSF ECCS-2154650.}
\thanks{$^{1}$Yujie Tang is with the Department of Industrial Engineering \& Management, Peking University, Beijing, 100871, China. 
        {yujietang@pku.edu.cn}.}%
\thanks{$^{2}$Yang Zheng is with the Department of Electrical and Computer Engineering, University of California San Diego, La Jolla, CA 92093, USA
        {zhengy@eng.ucsd.edu}.}%
}

\begin{document}

\maketitle
\thispagestyle{empty}
\pagestyle{empty}

%%%%%%%%%%%%%%%%%%%%%%%%%%%%%%%%%%%%%%%%%%%%%%%%%%%%%%%%%%%%%%%%%%%%%%%%%%%%%%%%
\begin{abstract}

Direct policy search has achieved great empirical success in reinforcement learning. Recently, there has been increasing interest in studying its theoretical properties for continuous control, and fruitful results have been established for linear quadratic regulator (LQR) and linear quadratic Gaussian (LQG) control that are smooth and nonconvex. In this paper, we consider the standard $\mathcal{H}_\infty$ robust control for output feedback systems and investigate the global optimality of direct policy search. Unlike LQR or LQG, the $\mathcal{H}_\infty$ cost function is nonsmooth in the policy space. Despite the lack of smoothness and convexity, our main result shows that for a class of \textit{non-degenerated} stabilizing controllers, all Clarke stationary points of $\mathcal{H}_\infty$ robust control 
are globally optimal and there is no spurious local minimum. Our proof technique is motivated by the idea of differentiable convex liftings (\texttt{DCL}), and we extend \texttt{DCL} to analyze the nonsmooth and nonconvex $\mathcal{H}_\infty$ robust control via convex reformulation. Our result sheds some light on the analysis of direct policy search for solving nonsmooth and nonconvex robust control problems.         

\end{abstract}

%%%%%%%%%%%%%%%%%%%%%%%%%%%%%%%%%%%%%%%%%%%%%%%%%%%%%%%%%%%%%%%%%%%%%%%%%%%%%%%%
\section{Introduction}

Inspired by the empirical successes of reinforcement~learning, direct policy search techniques have recently received extensive interest in the field of control. Significant advances have been established in terms of understanding the theoretical properties of direct policy search on a range of benchmark control problems, including stabilization~\cite{perdomo2021stabilizing}, linear quadratic regulator (LQR) \cite{fazel2018global,malik2019derivative,mohammadi2021convergence}, linear risk-sensitive control \cite{zhang2021policy},~linear quadratic Gaussian (LQG) control \cite{zheng2021analysis,zheng2022escaping,duan2022optimization}, dynamic filtering~\cite{umenberger2022globally}, and linear distributed control \cite{furieri2020learning,li2022distributed}; see \cite{hu2022towards} for a recent survey.   

All these control problems are known to be nonconvex in the policy space. One typical approach to deal with the nonconvexity in classical control theory is to reparameterize the problem into a convex form, e.g. via a suitable change of variables \cite{scherer1997multiobjective,zheng2021equivalence}, for which efficient algorithms exist \cite{zhou1996robust,boyd1994linear}. The reparameterization often relies on the underlying system model explicitly and is thus a model-based design. On the other hand, despite the nonconvexity, a series of recent findings have revealed favorable optimization landscape properties in many benchmark control problems. For example, global convergence of model-free policy gradient methods has been established for both discrete-time~\cite{fazel2018global} and continuous-time LQR~\cite{mohammadi2021convergence} thanks to the \textit{gradient dominance} property of the cost functions; the LQG cost function has no spurious stationary points that correspond to controllable and observable controllers~\cite{zheng2021analysis}. Beyond LQR and LQG, global or local convergence results of direct policy search have also been established for linear risk-sensitive control \cite{zhang2021policy} and distributed control problems \cite{furieri2020learning,li2022distributed}.    

For all the aforementioned benchmark control problems, the cost functions are typically linear and quadratic~in~terms of the system trajectories. Their policy optimization formulations thus admit a smooth cost function over the feasible region. This basic fact underpins most of the~existing landscape results and convergence guarantees~\cite{perdomo2021stabilizing,fazel2018global,malik2019derivative,mohammadi2021convergence,zhang2021policy,zheng2021analysis,zheng2022escaping,umenberger2022globally,furieri2020learning,li2022distributed,duan2022optimization}.
In addition to linear quadratic (LQ) control, another fundamental control paradigm, known as \emph{robust control}, addresses the worst-case performance against uncertainties~\cite{zhou1996robust}. In this case, the performance measure is the $\mathcal{H}_\infty$ norm of a certain closed-loop transfer function. It is known that the closed-loop $\mathcal{H}_\infty$ norm is not always differentiable in the policy space~\cite{apkarian2006nonsmooth}, which requires techniques from nonsmooth analysis to investigate the behavior of direct policy search. Indeed, a large class of nonsmooth policy search algorithms has been proposed~\cite{apkarian2006nonsmooth,apkarian2006nonsmooth2,apkarian2008mixed}, but these studies do not address the global optimality of direct policy search. The most relevant studies are \cite{guo2022global,hu2022connectivity}. In particular, the work \cite{guo2022global} has established a global convergence result of direct policy search for state-feedback $\mathcal{H}_\infty$ control. The work \cite{hu2022connectivity} considers output-feedback $\mathcal{H}_\infty$ dynamic control and reveals that there always exists a continuous path connecting any initial stabilizing controller to a global optimal controller.    

In this work, we focus on the policy optimization perspective for output-feedback $\mathcal{H}_\infty$ dynamic control. The $\mathcal{H}_\infty$ policy optimization formulation is known to be nonsmooth and nonconvex \cite{apkarian2006nonsmooth}. Our main result reveals that Clarke stationary points corresponding to \textit{non-degenerate dynamic controllers} are globally optimal (\Cref{theorem:main}); the precise definition of {non-degenerate dynamic controllers} will be discussed in Section~\ref{section:results}.
Our analysis employs tools from the classical convex reformulation of $\mathcal{H}_\infty$ control~\cite{scherer1997multiobjective,masubuchi1998lmi} and is also inspired by the recent works~\cite{umenberger2022globally,guo2022global}. Especially, we extend and tailor the technique of differentiable convex liftings (\texttt{DCL}) in~\cite{umenberger2022globally} to the nonsmooth setting. Our results and analysis shed light on policy optimization methods for solving nonsmooth and nonconvex robust control problems.

The rest of the paper is structured as follows. We present preliminaries and policy optimization formulation of $\mathcal{H}_\infty$ control in \Cref{section:preliminaries}. Our main results on global optimality are presented in \Cref{section:results}, and the technical proofs are provided in \Cref{section:proofs}. We show numerical experiments in \Cref{sec:examples} and conclude the paper in \Cref{section:conclusion}.

\noindent {\textit{Notations.}} We use $\mathbb{S}_{++}^n$ to denote the set of $n\times n$ real symmetric and positive definite matrices, and use $\mathrm{GL}_{n}$ to denote the set of $n\times n$ real invertible matrices. For two real symmetric matrices $A,B$, we use $A\succeq B$ or $B\preceq A$ (resp. $A\succ B$ or $B\prec A$) to mean that the matrix $A-B$ is positive semidefinite (reps. positive definite).

\section{Preliminaries and Problem Formulation} \label{section:preliminaries}

\subsection{Formulation of $\mathcal{H}_\infty$ Control}
We consider a continuous-time linear time-invariant (LTI) system of the form
\begin{equation}\label{eq:plant}
\begin{aligned}
\dot{x}(t) =\ & Ax(t) + B_1w(t) + B_2 u(t), \\
z(t) =\ & C_1 x(t) + D_{11}w(t) + D_{12} u(t), \\
y(t) =\ & C_2x(t) + D_{21} w(t),
\end{aligned}
\end{equation}
where $x(t)\in\mathbb{R}^{n_x}$ is the state of the plant, $w(t)\in\mathbb{R}^{n_w}$ represents exogenous disturbance, $u(t)\in\mathbb{R}^{n_u}$ is the control input, $z(t)\in\mathbb{R}^{n_z}$ represents the regulated performance signal, and $y(t)\in\mathbb{R}^{n_y}$ is the measured output. We make the following standard assumption. 
\begin{assumption}
$(A,B_1)$ and $(A,B_2)$ are controllable, and  $(C_1,A)$ and $(C_2,A)$ are observable.
\end{assumption}
 A typical control task is to synthesize a feedback controller (policy) that maps the output $y(t)$ to the control input $u(t)$, which stabilizes the plant and minimizes a certain performance metric. When we only have access to the output signal $y(t)$, a static feedback policy is typically not sufficient to ensure good closed-loop performance \cite{zhou1996robust}. We consider a full-order dynamic feedback policy of the form
\begin{equation} \label{eq:dynamic-controller}
    \begin{aligned}
\dot{\xi}(t) =\ &
A_\mK \xi(t) + B_\mK y(t), \\
u(t) =\ &
C_\mK \xi(t) + D_\mK y(t),
\end{aligned}
\end{equation}
where $\xi(t) \in\mathbb{R}^{n_x}$ represents the internal state of the dynamic controller. Then it is not difficult to see that the closed-loop transfer function from the exogenous disturbance $w$ to the output $z$ is given by
\begin{equation} \label{eq:transfer-function-Tzw}
    \mathbf{T}_{zw}(s) = C_{\mathrm{cl}}(\mK)
\left(sI - A_{\mathrm{cl}}(\mK)
\right)^{-1}
B_{\mathrm{cl}}(\mK)
+D_{\mathrm{cl}}(\mK),
\end{equation}
where we denote
\begin{align*}
A_{\mathrm{cl}}(\mK)
\coloneqq\  &
\begin{bmatrix}
A + B_2D_\mK C_2 & B_2C_\mK \\ B_\mK C_2 & A_\mK
\end{bmatrix}, \\
B_{\mathrm{cl}}(\mK)
\coloneqq\  & \begin{bmatrix}
B_1 + B_2D_\mK D_{21} \\
B_\mK D_{21}
\end{bmatrix}, \\
C_{\mathrm{cl}}(\mK)
\coloneqq\  &
\begin{bmatrix}
C_1 + D_{12}D_\mK C_2 & D_{12}C_\mK
\end{bmatrix}, \\
D_\mathrm{cl}(\mK) \coloneqq\ &
D_{11} + D_{12}D_\mK D_{21}.
\end{align*}
The goal of $\mathcal{H}_\infty$ control \cite{zhou1996robust} is to find a controller~\eqref{eq:dynamic-controller}~to~minimize the $\mathcal{H}_\infty$ norm of the transfer function $\mathbf{T}_{zw}$ defined~by
\begin{equation} \label{eq:Hinf-Tzw} 
\|\mathbf{T}_{zw}\|_\infty \coloneqq \sup_{\omega\in\mathbb{R}}\sigma_{\max}(\mathbf{T}_{zw}(j\omega)),
\end{equation}
 where $\sigma_{\max}(\cdot)$ denotes the largest singular value.

$\mathcal{H}_\infty$ control is a classical problem in control theory. Different approaches have been developed to find a suboptimal $\mathcal{H}_\infty$ controller, including the Riccati-equation-based approach \cite{doyle1988state} and the linear matrix inequality (LMI)-based approach \cite{gahinet1994linear}. Unlike $\mathcal{H}_2$ optimal control, the problem of finding an \emph{optimal} $\mathcal{H}_\infty$ controller does not have a closed-form solution\footnote{Given any $\gamma$ that is greater than the infimal value of $\|\mathbf{T}_{zw}\|_\infty$, one can find a controller satisfying $\|\mathbf{T}_{zw}\|_\infty \leq \gamma$ by employing the Riccati-equation-based approach \cite{doyle1988state}. However, this controller in general is a \emph{suboptimal} controller and does not achieve the optimal value of $\|\mathbf{T}_{zw}\|_\infty$.
}.

\subsection{Problem Statement: Nonsmooth Policy Optimization}
In this paper, we investigate the perspective of policy optimization that directly searches over the (parameterized) controller/policy space. We parameterize the dynamic controller by its associated system matrices as
$$
\mathsf{K} = \begin{bmatrix}
D_\mK & C_\mK  \\
B_\mK & A_\mK
\end{bmatrix} \in \mathbb{R}^{(n_u + n_x) \times (n_y + n_x)},
$$
where we fix the dimension of $\xi(t)$ to be $n_x$. Define
\begin{equation} \label{eq:stabilizing-controller}
\mathcal{C}
\coloneqq\ 
\left\{
\mK: A_{\mathrm{cl}}(\mK)\text{ is stable}
\right\}.
\end{equation}
The closed-loop system is internally stable if and only if $\mK\in\mathcal{C}$ \cite{zhou1996robust}. Therefore, any internally stabilizing controller can be parameterized by $\mK\in\mathcal{C}$. The $\mathcal{H}_\infty$ control problem can then be reformulated as policy optimization of the form
\begin{equation}\label{eq:H_infty_policy_opt}
\min_{\mK}\ \ J(\mK) \quad\text{subject to}\ \ \mK\in\mathcal{C},
\end{equation}
where $J(\mK)$ denotes the $\mathcal{H}_\infty$ norm $\|\mathbf{T}_{zw}\|_\infty$ for each $\mK\in\mathcal{C}$. 

The idea of direct policy search is to start from an initial policy $\mK_0 \in \mathcal{C}$ and conduct the iteration $\mK_{t+1} = \mK_t + \alpha_t F_t, t\geq 0$, where $\alpha_t >0$ is a step size and $F_t$ is a search direction, such that the $\mathcal{H}_\infty$ cost $J(\mK_t)$ is gradually improved. However, the $\mathcal{H}_\infty$ cost function $J(\mK)$ in \eqref{eq:Hinf-Tzw} is known to be nonconvex and also nonsmooth with two possible sources of non-smoothness: One from taking the largest singular of complex matrices, and the other from maximization over all the frequencies $\omega \in \mathbb{R}$.
For the nonconvex and nonsmooth $\mathcal{H}_\infty$ policy optimization problem~\eqref{eq:H_infty_policy_opt}, it is unclear yet how to certify the global optimality of direct policy search methods. This motivates the main question of our work:  

\vspace{-1mm}

\begin{center}
    \textit{Can we characterize a class of stationary points that are globally optimal to \eqref{eq:H_infty_policy_opt} despite the non-convexity and non-smoothness?}  
\end{center}

\vspace{-1mm}

 This paper provides a positive answer to this question. Our analysis employs the classical convex reformulation \cite{scherer1997multiobjective,masubuchi1998lmi}, and is also motivated by the idea of differentiable convex liftings (\texttt{DCL}) for dynamic filtering~\cite{umenberger2022globally} and the study on state-feedback $\mathcal{H}_\infty$ control \cite{guo2022global}. We make non-trivial extensions to dynamic output feedback $\mathcal{H}_\infty$ control \eqref{eq:H_infty_policy_opt}.

\begin{remark}
    Direct policy search has been used in earlier studies~\cite{apkarian2006nonsmooth,apkarian2006nonsmooth2,apkarian2008mixed} to address $\mathcal{H}_\infty$ controller synthesis, but no optimality~guarantees are given. It has regained increasing attention due to recent advances in policy optimization for learning-based control~\cite{hu2022connectivity,zheng2021analysis,umenberger2022globally,guo2022global}. In particular, the recent work~\cite{guo2022global} has established a global optimality guarantee of policy search for state-feedback $\mathcal{H}_\infty$ control. Apart from better scalability compared to the classical approaches~\cite{apkarian2006nonsmooth,apkarian2006nonsmooth2}, the direct policy search approach appears more amenable to model-free control when a model of the plant is lacking.
\end{remark}

\subsection{Clarke Subdifferential}

We review the notion of Clarke subdifferential \cite{clarke1990optimization} that enables the analysis of a large class of nonsmooth functions. 
Let $f(x): C\rightarrow\mathbb{R}$ be a function defined on an open subset $C\subseteq\mathbb{R}^n$. We say that $f$ is \textit{locally Lipschitz near $x\in C$}, if there exists $\epsilon>0$ and $L>0$ such that for any $y_1,y_2\in C$ satisfying $\|y_1-x\|<\epsilon$ and $\|y_2-x\|<\epsilon$, we have $|f(y_1) - f(y_2)| \leq L \|y_1 - y_2\|$. The function $f$ is said to be \textit{locally Lipschitz over $C$} if it is locally Lipschitz near any $x\in C$. The Rademacher theorem \cite[Theorem 3.2]{evans2015measure} guarantees that a locally Lipschitz function is differentiable almost everywhere in the domain.

Let $f(x)$ be locally Lipschitz over $C$. We define its Clarke directional derivative at $x\in C$ in the direction $v\in\mathbb{R}^n$ by
\[
f^\circ(x;v)
=\limsup_{x'\rightarrow x,t\downarrow 0}
\frac{f(x'+tv)-f(x')}{t}.
\]
The local Lipschitz continuity of $f$ guarantees that $f^\circ(x;v)$ is finite for all $x\in C$ and $v\in\mathbb{R}^n$. It can be shown that for any fixed $x\in C$, $f^\circ(x;\cdot)$ is a convex function and satisfies $f^\circ(x;\lambda v)=\lambda f^\circ(x;v)$ for any $\lambda>0$. We then define the Clarke subdifferential of $f$ at $x\in C$ as the set
\[
\partial f(x) \coloneqq
\left\{g\in\mathbb{R}^n: f^\circ(x;v)\geq\langle g,v\rangle\text{ for all }v\in\mathbb{R}^n\right\},
\]
which is nonempty for any $x\in C$.
It is shown that $f^\circ(x,\cdot)$ is the support function of $\partial f(x)$  \cite[Proposition 2.1.2]{clarke1990optimization}:
\begin{equation} \label{eq:support-function}
f^\circ(x,v) = \max_{g \in \partial f(x)}\; \langle g,v\rangle. 
\end{equation}
Moreover, the following equality holds~\cite[Theorem 2.5.1]{clarke1990optimization}:
\begin{equation*}
    \partial f(x) = \operatorname{conv}\left\{\lim_{ x_k \to x} \nabla f(x_k) \,\Big|\, \nabla f(x_k) \text{ exists}, x_k \in C\right\},
\end{equation*}
where $\mathrm{conv}$ denotes the convex hull of a set.

We call $x \in C$ a \textit{Clarke stationary point} if $0 \in \partial f(x)$. The following result relates local minima and local maxima with Clarke stationary points \cite[Proposition 2.3.2]{clarke1990optimization}.
\begin{lemma}\label{lemma:clarke_stationary}
Let $f(x)$ be locally Lipschitz over $C$. If $x\in C$ is a local minimum or maximum of $f(x)$, then $x$ is a Clarke stationary point, i.e., $0 \in \partial f(x)$. 
\end{lemma}
Note that the converse of \Cref{lemma:clarke_stationary} does not hold in general.   
The function $f$ is called \textit{subdifferntial regular}, if for any $x \in C$, the ordinary directional derivative exists and coincides with the Clarke directional~derivative for all directions, i.e.,
\[
\lim_{t\downarrow 0}
\frac{f(x+tv)-f(x)}{t}
=
f^\circ(x,v),\quad\forall v \in \mathbb{R}^n,x\in C.
\]
We denote the ordinary directional derivative by $f'(x;v)$ whenever it exists. 
From~\eqref{eq:support-function}, the following result is clear. 
\begin{lemma}\label{lemma:regular_Clarke_stationary}
    Suppose that $f$ is subdifferential regular. If $x$ is a Clarke stationary point, then $f'(x,v) \geq 0$ for all $v\in \mathbb{R}^n$. 
\end{lemma}

\section{Main Results} \label{section:results}

In this section, we first summarize some useful properties of the $\mathcal{H}_\infty$ cost function $J(\mK)$ in Section \ref{subsection:basic-properties}. Our main technical result is then presented in Section \ref{subsection:main-results}, which characterizes a large class of stationary points that are globally optimal. The rest of this section presents its proof.

\subsection{Basic Properties of the $\mathcal{H}_\infty$ Cost} \label{subsection:basic-properties}

We first review a fact for $\mathcal{H}_\infty$ policy optimization.
\begin{lemma}
    The set of internally stabilizing controllers~$\mathcal{C}$ is nonconvex, potentially disconnected, but has at most two connected components. The cost function $J(\mK)$ is nonconvex. 
\end{lemma}
We refer to~\cite{hu2022connectivity,zheng2021analysis} for relevant proofs and examples. The following result is also known in the community.

\begin{proposition}[{\!\!{\cite[Proposition 3.1]{apkarian2006nonsmooth2}}}]\label{lemma:H_inf_some_local_Lipschitz}
For the $\mathcal{H}_\infty$ policy optimization problem~\eqref{eq:H_infty_policy_opt}, the following statements hold.
\begin{enumerate}[leftmargin=12pt, labelwidth=8pt,labelsep=4pt]
\item $J(\mK)$ is locally Lipschitz over $\mK\in\mathcal{C}$.

\item $J(\mK)$ is subdifferentially regular.
\end{enumerate}
\end{proposition}
 The proof idea in {\cite[Proposition 3.1]{apkarian2006nonsmooth2}} is to view $J(\mK)$ as a composition of a convex mapping $\|\cdot\|_\infty$ and the mapping $\mK\mapsto \mathbf{T}_{zw}$ that is continuously differentiable over $\mK$. Then, the subdifferential regularity of $J(\mK)$ follows from \cite{clarke1990optimization}. We provide some missing details in Appendix A.

\Cref{lemma:H_inf_some_local_Lipschitz} justifies that $J(\mK)$ is Clarke subdifferentiable. It is now clear from \Cref{lemma:clarke_stationary} that if a dynamic controller $\mK \in \mathcal{C}$ is a local minimum of $J(\mK)$, then $\mK$ is a Clarke stationary point. Our main goal is to establish a class of Clarke stationary points that are globally optimal to~\eqref{eq:H_infty_policy_opt}. In our analysis, the following bounded-real lemma will be fundamental.

\begin{lemma}[Bounded real lemma]
\label{lemma:bounded_real}
Let $A$ be stable, and consider a transfer function
$
G(s) = C(sI-A)^{-1}B + D.
$
Let $\gamma>0$ be arbitrary. The following statements hold. 

\begin{enumerate}[leftmargin=12pt,labelwidth=8pt,labelsep=4pt]
\item (Strict version, \cite[Lemma 7.3]{dullerud2013course}) $\|G\|_{\infty}< \gamma$ if and only if there exists $P\succ 0$ such that
\begin{equation} \label{eq:strict-hinf}
\begin{bmatrix}
A^\tr P + P A & PB & C^\tr \\
B^\tr P & -\gamma I & D^\tr \\
C & D & -\gamma I
\end{bmatrix}\prec 0.
\end{equation} 

\item (Nonstrict version,
\cite[Section 2.7.3]{boyd1994linear}) $\|G\|_{\infty}\leq \gamma$ if there exists $P\succ 0$ such that
\begin{equation} \label{eq:non-strict-hinf}
    \begin{bmatrix}
A^\tr P + P A & PB & C^\tr \\
B^\tr P & -\gamma I & D^\tr \\
C & D & -\gamma I
\end{bmatrix}\preceq 0.
\end{equation}
The converse holds if $(A,B,C)$ is controllable and observable.
\end{enumerate}
\end{lemma}

\subsection{Main Technical Results} \label{subsection:main-results}
To state our main results, we introduce a special class of controllers in $\mathcal{C}$, which we will call \emph{non-degenerate} stabilizing controllers below. Specifically, we define
\begin{equation} \label{eq:Set-Snd}
\begin{aligned}
\mathcal{S}_{\mathrm{nd}}\!
\coloneqq\!
\Bigg\{\!(\mK,P,\gamma):&
\, P=\begin{bmatrix}
P_{11} & P_{12} \\
P_{12}^\tr & P_{22}
\end{bmatrix}
\in\mathbb{S}_{++}^{2n_x},\ \\
&\quad  P_{12}\in\mathrm{GL}_{n_x}, \;
\mathscr{N}\!(\mK,P,\gamma)
\!\preceq 0
\Bigg\}.
\end{aligned}
\end{equation}
where we denote
\begingroup
    \setlength\arraycolsep{2pt}
\def\arraystretch{0.9}
$$
\mathscr{N}\!(\mK,P,\gamma)
\!:=\!\begin{bmatrix}
A_{\mathrm{cl}}(\mK)^\tr P \!+\! P A_{\mathrm{cl}}(\mK) & PB_{\mathrm{cl}}(\mK) & C_{\mathrm{cl}}(\mK)^\tr \\
B_{\mathrm{cl}}(\mK)^\tr P & -\gamma I & D_{\mathrm{cl}}(\mK)^\tr \\
C_{\mathrm{cl}}(\mK) & D_{\mathrm{cl}}(\mK) & -\gamma I
\end{bmatrix}.
$$
\endgroup
It is clear that for any triple $(\mK,P,\gamma) \in \mathcal{S}_{\mathrm{nd}}$, we have $\mK \in \mathcal{C}$ and $J(\mK) \leq \gamma$ by the non-strict version of \Cref{lemma:bounded_real}. Note that when defining $\mathcal{S}_{\mathrm{nd}}$ in \eqref{eq:Set-Snd}, we require the off-diagonal block $P_{12}$ to have full rank, which will be explained in Remark~\ref{remark:invertibility_P12}. We further define
\begin{equation} \label{eq:Set-Cnd}
\begin{aligned} 
\mathcal{C}_{\mathrm{nd}}
\coloneqq
\left\{\mK\in\mathcal{C}:\exists P
\text{ such that }(\mK,P,J(\mK))\in\mathcal{S}_{\mathrm{nd}}\right\}.
\end{aligned} 
\end{equation}
Controllers in $\mathcal{C}_{\mathrm{nd}}$ will be called \emph{non-degenerate} stabilizing controllers, since each controller in $\mathcal{C}_{\mathrm{nd}}$ admits a $P$ with a non-degenerate off-diagonal block $P_{12}$ to certify the associated $\mathcal{H}_\infty$ cost $J(\mK)$ in \eqref{eq:non-strict-hinf}. 

By definition, we have $\mathcal{C}_{\mathrm{nd}} \subseteq \mathcal{C}$. We conjecture that non-degenerate stabilizing controllers are ``generic'' in the sense that the complement set $\mathcal{C}\backslash\mathcal{C}_{\mathrm{nd}}$ has measure zero. A rigorous proof of this conjecture seems challenging and is still ongoing work.
In \Cref{sec:examples}, we shall provide some numerical evidence of this conjecture.

{
\begin{remark}[Invertibility of $P_{12}$]\label{remark:invertibility_P12}
    In \eqref{eq:Set-Snd}, we require the off-diagonal block $P_{12}$ to have full rank. This requirement on $P_{12}$ is essential in deriving the convex reformulation of $\mathcal{H}_2$ or $\mathcal{H}_\infty$ control proposed in \cite{scherer1997multiobjective,masubuchi1998lmi}. On the other hand, when only strict LMIs (e.g., \eqref{eq:strict-hinf}) are imposed, we can slightly perturb $P$ to ensure that $P_{12}$ has full rank without violating the strict LMIs, which is a trick that has been employed in \cite{scherer1997multiobjective,masubuchi1998lmi} as well as some recent studies \cite{zheng2021analysis,hu2022connectivity}. But in this paper, we aim to directly analyze the $\mathcal{H}_\infty$ cost function $J(\mK)$ instead of its upper bound, meaning that our subsequent results and proofs need to use the non-strict LMI~\eqref{eq:non-strict-hinf}. Therefore, we need to explicitly require the off-diagonal block $P_{12}$ in \eqref{eq:Set-Snd} to be invertible. Similar requirements appear in the setting of dynamic filtering in \cite{umenberger2022globally}, which were called \textit{informativity} by the authors. 
\end{remark}
}

We are now ready to state our main technical result.
\begin{theorem}\label{theorem:main}
    Given any non-degenerate stabilizing controller $\mK\in\mathcal{C}_{\mathrm{nd}}$,
    if $\mK$ is a Clarke stationary point, i.e., $0 \in \partial J(\mK)$, then it is a global minimum of $J(\mK)$ over $\mathcal{C}$.
\end{theorem}

This result also highlights that there are no spurious local minima in the set of non-degenerate stabilizing controllers $\mK\in\mathcal{C}_{\mathrm{nd}}$. The following corollary is immediate.

\begin{corollary}
For the $\mathcal{H}_\infty$ policy optimization problem \eqref{eq:H_infty_policy_opt}, we have 
\begin{itemize}[leftmargin=12pt, labelwidth=6pt,labelsep=6pt]
    \item Any local minimum of $J(\mK)$ in $\mathcal{C}_{\mathrm{nd}}$ is a global minimum. 
    \item There exists no local maximum of $J(\mK)$ in $\mathcal{C}_{\mathrm{nd}}$.
\end{itemize}
\end{corollary}

\begin{remark}
It is known that the feasible region of~\eqref{eq:H_infty_policy_opt} has at most two connected components \cite{zheng2021analysis}.
Moreover, \cite{hu2022connectivity} has also shown that there always exists a continuous path from any initial point $\mK_0 \in \mathcal{C}$ to a global minimum. Thus it makes no difference to search over either connected component in $\mathcal{C}$ when solving \eqref{eq:H_infty_policy_opt} via direct policy search. Our result in \Cref{theorem:main} has further provided a global optimality certificate for \eqref{eq:H_infty_policy_opt}, showing positive news for direct policy search methods. Note that any stationary points corresponding to \textit{controllable and observable controllers} in $\mathcal{H}_2$ control are globally optimal \cite[Theorem 4.3]{zheng2021analysis}.  \Cref{theorem:main} can thus be viewed as the counterpart in output-feedback $\mathcal{H}_\infty$ control.
\end{remark}

The proof of \Cref{theorem:main} was inspired by the idea of differentiable convex liftings (\texttt{DCL}) for output estimation~\cite{umenberger2022globally} and the recent analysis on state-feedback $\mathcal{H}_\infty$ control \cite{guo2022global}. In this paper, we make non-trivial extensions of the \texttt{DCL} analysis to the nonsmooth output feedback $\mathcal{H}_\infty$  control problem.
The following subsection gives the proof of \Cref{theorem:main}.

\subsection{Proof of \Cref{theorem:main}} \label{subsection:proof-Theorem-1}

\begin{figure*}[t]
  \hrulefill
  \begin{equation}\label{eq:bigM_def}
\begin{aligned}
& \mathscr{M}(X,Y,M,H,F,G,\gamma) \\
\coloneqq &
\begin{bmatrix}
AX \!+\! B_2F \!+\! (AX \!+\! B_2F)^\tr & M^\tr \!+\! A \!+\! B_2GC_2 & B_1 \!+\! B_2GD_{21} & (C_1 X \!+\! D_{12}F)^\tr \\
M \!+\! (A \!+\! B_2GC_2)^\tr & YA \!+\! HC_2 \!+\! (YA \!+\! HC_2)^\tr & YB_1 \!+\! HD_{21} & (C_1 \!+\! D_{12}GC_2)^\tr \\
(B_1 \!+\! B_2GD_{21})^\tr & (YB_1 \!+\! HD_{21})^\tr & -\gamma I & (D_{11} \!+\! D_{12}GD_{21})^\tr \\
C_1 X \!+\! D_{12}F & C_1 \!+\! D_{12}GC_2 & D_{11} \!+\! D_{12}GD_{21} & -\gamma I
\end{bmatrix},
\end{aligned}
  \end{equation}
  \hrulefill
  \vspace{-10pt}
\end{figure*}

We first introduce some auxiliary quantities.
Given a set of matrices $X\!\in \mathbb{S}^{n_x},Y \!\in \mathbb{S}^{n_x},
M \!\in\! \mathbb{R}^{n_x\times n_x},
H \!\in\! \mathbb{R}^{n_x\times n_y},
F \!\in\! \mathbb{R}^{n_u\times n_x},
G \!\in\! \mathbb{R}^{n_u\times n_y},
\gamma \!\in\! \mathbb{R}$, we define an affine function $\mathscr{M}(X,Y,M,H,F,G,\gamma)$ by \eqref{eq:bigM_def}, and then define a convex set
\begin{align}  \nonumber
\mathcal{F}\!= \! &
\Bigg\{ \! (X, Y, M, H, F, G,\gamma)\!
:
X,Y \!\in \mathbb{S}^{n_x},
M \!\in\! \mathbb{R}^{n_x\times n_x}, \\ \nonumber
& \quad
H \!\in\! \mathbb{R}^{n_x\times n_y},
F \!\in\! \mathbb{R}^{n_u\times n_x},
G \!\in\! \mathbb{R}^{n_u\times n_y},
\gamma \!\in\! \mathbb{R},\\
& \quad
\begin{bmatrix}
X & I \\
I & Y
\end{bmatrix}\succ 0,
\mathscr{M}(X,Y,M,H,F,G,\gamma)\preceq 0
\Bigg\},
\end{align}
and an extended set 
\begin{equation}
\mathcal{G} = \mathrm{GL}_{n_x}\times \mathcal{F}.
\end{equation}
We note that the LMI $\mathscr{M}(X,Y,M,H,F,G,\gamma)\preceq 0$ resembles the structure in the non-strict bounded real lemma \eqref{eq:non-strict-hinf} as well as the LMI in $\mathcal{S}_{\textrm{nd}}$ \cref{eq:Set-Snd}. Indeed, based on a non-trivial change of variables in \cite{scherer1997multiobjective} that reformulates an output feedback $\mathcal{H}_\infty$ control problem into a set of LMIs, we can build a smooth bijection between the set $\mathcal{S}_{\mathrm{nd}}$ and set $\mathcal{G}$.

In particular, for each $(\mK,P,\gamma)\in\mathcal{S}_{\mathrm{nd}}$, we define the mapping $\Phi(\mK,P,\gamma)$ by
\begin{align}
\Phi(\mK,P,\gamma) \!= 
\big(P_{12},(P^{-1})_{11},\!
P_{11},
\Phi_M,\! \Phi_H,\!
\Phi_F,\!
D_{\mK},\!
\gamma\big),
\end{align}
where
\begin{align*}
\Phi_M \coloneqq\ &
P_{12}B_\mK C_2 (P^{-1})_{11}
+ P_{11} B_2 C_\mK (P^{-1})_{21} \\
& + P_{11}(A+B_2 D_\mK C_2)(P^{-1})_{11} 
 + P_{12}A_\mK (P^{-1})_{21}, \\
\Phi_{H} \coloneqq\ &
P_{11}B_2 D_\mK + P_{12} B_\mK, \\
\Phi_{F} \coloneqq\ &
D_\mK C_2 (P^{-1})_{11} + C_\mK (P^{-1})_{21},
\end{align*}
and $(P^{-1})_{12}$, for instance, denotes the 
$n_x\times n_x$ submatrix of $P^{-1}$ corresponding to the first $n_x$ row and last $n_x$ columns. 

We have the following result that shows the connection between the sets $\mathcal{S}_{\mathrm{nd}}$, $\mathcal{G}$, and the mapping $\Phi$.

\begin{proposition} \label{proposition:bijection}
$\Phi$ is a diffeomorphism from $\mathcal{S}_{\mathrm{nd}}$ to $\mathcal{G}$, i.e., $\Phi$ is indefinitely differentiable and invertible, and $\Phi^{-1}$ is also indefinitely differentiable.
\end{proposition}

The proof of this proposition is mostly based on direct constructions which are motivated by the change of variables in \cite{scherer1997multiobjective}. We first notice that each element of $\Phi$ is a rational function over the domain $\mathcal{S}_{\mathrm{nd}}$, and thus $\Phi$ is real analytic. By direct verification, we can show that $\Phi$ maps $\mathcal{S}_{\mathrm{nd}}$ into $\mathcal{G}$. Further, we can explicitly construct the inverse mapping of $\Phi$, which is also real analytic. This proves that $\Phi$ is a diffeomorphism from $\mathcal{S}_{\mathrm{nd}}$ to $\mathcal{G}$. The detailed steps are provided in \Cref{subsection:proof-bijection}.  

After establishing the connection between $\mathcal{S}_{\mathrm{nd}}$ and $\mathcal{G}$ via the mapping $\Phi$, we can further derive the following two technical results. Their proofs are inspired by the recently proposed framework of \texttt{DCL}, but we extend and tailor the relevant techniques to the nonsmooth $\mathcal{H}_\infty$ control setting. 
The details are technically involved, and we postpone them to \Cref{subsection:proof-before-main,subsection:proposition-main}. 

\begin{proposition}
\label{lemma:intermediate_before_min}
Let $(\mK, P, \gamma)\in \mathcal{S}_\mathrm{nd}$ be arbitrary, and suppose there exists $(\mK',P',\gamma')\in\mathcal{S}_{\mathrm{nd}}$ such that $\gamma>\gamma'$. Then there exists a $C^\infty$ curve $\psi:[0,\delta)\rightarrow \mathcal{S}_{\mathrm{nd}}$ satisfying $\psi(0)=(\mK,P,\gamma)$
such that
\[
\lim_{t\downarrow 0}
\frac{\pi_\gamma(\psi(t))-\pi_\gamma(\psi(0))}{t}< 0,
\]
where $\pi_\gamma:\mathcal{S}_{\mathrm{nd}}\rightarrow\mathbb{R}$ denotes the canonical projection $\pi_\gamma(\mK,P,\gamma)=\gamma$.
\end{proposition}

\begin{proposition}\label{proposition:main}
Let $\mK\in\mathcal{C}_{\mathrm{nd}}$ be arbitrary, and suppose there exists $\mK'\in\mathcal{C}$ such that $J(\mK)>J(\mK')$. Then there exists $\mathsf{V}\neq 0$ such that
$$
\lim_{t\downarrow 0} \frac{J(\mK+t\mathsf{V})-J(\mK)}{t}
<0,
$$
i.e., the ordinary directional derivative of $J$ at $\mK$ in the direction $\mathsf{V}$ is strictly negative.
\end{proposition}

The proof of \Cref{theorem:main} becomes immediate by combining \Cref{proposition:main} with~\Cref{lemma:H_inf_some_local_Lipschitz} and~\Cref{lemma:regular_Clarke_stationary}. Indeed, \Cref{lemma:H_inf_some_local_Lipschitz} confirms that $J(\mK)$ is subdifferentially regular, and then \Cref{lemma:regular_Clarke_stationary} states that for any Clarke stationary point $\mK$, we have $J'(\mK,\mathsf{V}) \geq 0$ for all directions $\mathsf{V}$. Now consider a Clarke stationary point $\mK \in \mathcal{C}_{\mathrm{nd}}$. If it is not a globally minimum, then there exists another controller $\mK'\in\mathcal{C}$ such that $J(\mK)>J(\mK')$. Then, \Cref{proposition:main} guarantees that $J'(\mK,\mathsf{V}) < 0$ for some direction $\mathsf{V}$, which contradicts to \Cref{lemma:regular_Clarke_stationary}. Therefore, a Clarke stationary point $\mK \in \mathcal{C}_{\mathrm{nd}}$ must be a global minimum of $J(\mK)$.

\section{Technical proofs} \label{section:proofs}

\subsection{Proof of \Cref{proposition:bijection}} \label{subsection:proof-bijection}
We first show that $\Phi$ maps $\mathcal{S}_{\mathrm{nd}}$ into $\mathcal{G}$. Let $(\mK,P,\gamma)\in\mathcal{S}_{\mathrm{nd}}$ be arbitrary, and denote
\begin{align*}
{\Xi=}\ & P_{12},\ \ X = (P^{-1})_{11},\ \ Y = P_{11}, \\
M =\ & \Phi_M(\mK,P,\gamma),\ \ 
H = \Phi_H(\mK,P,\gamma), \\
F =\ & \Phi_F(\mK,P,\gamma),\ \ 
G = D_\mK,
\end{align*}
i.e., $(\Xi,X,Y,M,H,F,G,\gamma)=\Phi(\mK,P,\gamma)$. Let
\[
T = \begin{bmatrix}
(P^{-1})_{11} & I \\ (P^{-1})_{21} & 0
\end{bmatrix}.
\]
The definition of $\mathcal{S}_{\mathrm{nd}}$ directly implies that $\Xi=P_{12}\in\mathrm{GL}_{n_x}$. Also, $P\succ 0$ implies that $P$ is invertible, and we can infer from $PP^{-1}=I$ that
\begin{equation}\label{eq:prop2_intermediate0}
PT=\begin{bmatrix}
I & P_{11} \\ 0 & P_{12}^\tr
\end{bmatrix}.
\end{equation}
$P_{12}\in\mathrm{GL}_{n_x}$ then implies that $PT$ is invertible, which further implies that $T$ is invertible. Consequently,
\begin{align*}
0\prec\ &
T^\tr PT =
(PT)^\tr T
=\begin{bmatrix}
I & 0 \\ P_{11} & P_{12}
\end{bmatrix}
\begin{bmatrix}
(P^{-1})_{11} & I \\ (P^{-1})_{21} & 0
\end{bmatrix} \\
=\ &
\begin{bmatrix}
(P^{-1})_{11} & I \\
P_{11}(P^{-1})_{11}
+P_{12} (P^{-1})_{21} & P_{11}
\end{bmatrix}
=\begin{bmatrix}
X & I \\ I & Y
\end{bmatrix}.
\end{align*}
To show that $\mathscr{M}(X,Y,M,H,F,G,\gamma)\preceq 0$, we note that
\begin{equation}\label{eq:prop2_intermediate1}
\begin{aligned}
& T^\tr PA_{\mathrm{cl}}(\mK)T \\
=\ &
\begin{bmatrix}
I & 0 \\ P_{11} & P_{12}
\end{bmatrix}
\begin{bmatrix}
A+B_2D_\mK C_2 & B_2 C_\mK \\
B_\mK C_2 & A_\mK
\end{bmatrix}
\begin{bmatrix}
(P^{-1})_{11} & I \\ (P^{-1})_{12}^\tr & 0
\end{bmatrix} \\
=\ &
\begin{bmatrix}
AX + B_2F & A+B_2GC_2 \\
M & YA + HC_2
\end{bmatrix}.
\end{aligned}
\end{equation}
Similarly, it can be  verified that
\begin{align}
T^\tr PB_{\mathrm{cl}}(\mK)
=\ &
\begin{bmatrix}
B_1 + B_2GD_{21}  \\
YB_1 + HD_{21}
\end{bmatrix}, \label{eq:prop2_intermediate2} \\
C_{\mathrm{cl}}(\mK)T
=\ &
\begin{bmatrix}
C_1X + D_{12}F & 
C_1 + D_{12}GC_2
\end{bmatrix}. \label{eq:prop2_intermediate3}
\end{align}
Summarizing these identities, we can show that
\begin{align*}
& \mathscr{M}(X,Y,M,H,F,G,\gamma) \\
=\ &
\begin{bmatrix}
T & 0 & 0 \\
0 & I & 0 \\
0 & 0 & I
\end{bmatrix}^\tr
\mathscr{N}(\mK,P,\gamma)
\begin{bmatrix}
T & 0 & 0 \\
0 & I & 0 \\
0 & 0 & I
\end{bmatrix}.
\end{align*}
Since $(\mK,P,\gamma)\in\mathcal{S}_{\mathrm{nd}}$ implies $\mathscr{N}(\mK,P,\gamma)\preceq 0$, we get $\mathscr{M}(X,Y,M,H,F,G,\gamma)\preceq 0$. We can now conclude that $(\Xi,X,Y,M,H,F,G,\gamma)=\Phi(\mK,P,\gamma)\in\mathcal{G}$, and further $\Phi(\mathcal{S}_{\mathrm{nd}})\subseteq\mathcal{G}$.

We then show that $\Phi$ is a bijection from $\mathcal{S}_{\mathrm{nd}}$ to $\mathcal{G}$. We construct a mapping $\Psi$ defined on $\mathcal{G}$ as follows: Let $\Xi\in\mathrm{GL}_{n_x}$ and $\mathsf{Z}=(X,Y,M,H,F,G,\gamma)\in\mathcal{F}$ be arbitrary, and denote $\Pi = -\Xi^{-1}(Y-X^{-1})X$. Since $\begin{bmatrix}
X & I \\ I & Y
\end{bmatrix}
\succ 0$, we see that $Y-X^{-1}\succ 0$, which implies that $\Pi$ is invertible. Furthermore, we have
$
\Xi^\tr(Y-X^{-1})^{-1}\Xi \succ 0
$
and
\begin{align*}
& Y-\Xi\left(\Xi^\tr(Y-X^{-1})^{-1}\Xi\right)^{-1}\Xi^\tr
= X^{-1}\succ 0,
\end{align*}
which imply that
$$
\begin{bmatrix}
Y & \Xi \\
\Xi^\tr & \Xi^\tr(Y-X^{-1})^{-1}\Xi
\end{bmatrix}\succ 0.
$$
Now, we let
\begin{align*}
\Psi_\mK(\Xi,\mathsf{Z})
= &
\begin{bmatrix}
I & 0 \\
YB_2 & \Xi
\end{bmatrix}^{-1}\!
\begin{bmatrix}
G & F \\ H & M \!-\! YAX
\end{bmatrix}\!
\begin{bmatrix}
I & C_2 X \\ 0 & \Pi
\end{bmatrix}^{-1},
\\
\Psi_P(\Xi,\mathsf{Z})
= &
\begin{bmatrix}
Y & \Xi \\
\Xi^\tr & \Xi^\tr(Y-X^{-1})^{-1}\Xi
\end{bmatrix},
\end{align*}
and
$$
\Psi(\Xi,\mathsf{Z}) = \left(\Psi_\mK(\Xi,\mathsf{Z}),
\Psi_P(\Xi,\mathsf{Z}),\gamma\right).
$$
By definition, $\Psi_P(\Xi,\mathsf{Z})$ is positive definite, and the $(1,2)$-block of $\Psi_P(\Xi,\mathsf{Z})$ is invertible. Moreover,
we can verify
\[
\Psi_P(\Xi,\mathsf{Z}) \begin{bmatrix}
X & I \\ \Pi & 0
\end{bmatrix}
=\begin{bmatrix}
I & Y \\ 0 & \Xi^\tr
\end{bmatrix}.
\]
This equality has the same form as~\eqref{eq:prop2_intermediate0}. Therefore, we can mimic the calculations in deriving~\eqref{eq:prop2_intermediate1} to~\eqref{eq:prop2_intermediate3} to show that
\begin{align*}
& \begin{bmatrix}
X & I \\
\Pi & 0 \\
 & & I \\
 & & & I
\end{bmatrix}^{-\tr}\mathscr{M}(\mathsf{Z})
\begin{bmatrix}
X & I \\
\Pi & 0 \\
 & & I \\
 & & & I
\end{bmatrix}^{-1}\\
=\ &
\mathscr{N}(\Psi_\mK(\Xi,\mathsf{Z},\gamma),\Psi_P(\Xi,\mathsf{Z},\gamma),\gamma),
\end{align*}
implying that $\mathscr{N}(\Psi_\mK(\Xi,\mathsf{Z},\gamma),\Psi_P(\Xi,\mathsf{Z},\gamma),\gamma)\preceq 0$ whenever $(\Xi,\mathsf{Z})\in\mathcal{G}$. Thus we can conclude that $\Psi$ is a mapping from $\mathcal{G}$ into $\mathcal{S}_{\mathrm{nd}}$. We can then compute the compositions of mappings $\Psi\circ\Phi$ and $\Phi\circ\Psi$ by tedious but straightforward calculations, which turn out to be the identity maps on $\mathcal{S}_{\mathrm{nd}}$ and $\mathcal{G}$, respectively. Therefore $\Phi$ is a bijection from $\mathcal{S}_{\mathrm{nd}}$ to $\mathcal{G}$ with $\Psi$ being its inverse.

Finally, note that $\Phi$ and $\Psi$ are both real analytic over their domains. Thus $\Phi$ is a diffeomorphism from $\mathcal{S}_{\mathrm{nd}}$ to $\mathcal{G}$.

\subsection{Proof of \Cref{lemma:intermediate_before_min}} \label{subsection:proof-before-main}
Let $(\Xi,\mathsf{Z})=\Phi(\mK,P,\gamma)$ and $(\Xi',\mathsf{Z}') = \Phi(\mK',P',\gamma')$ where $\Xi,\Xi'\in\mathrm{GL}_{n_x}$ and $\mathsf{Z},\mathsf{Z}'\in\mathcal{F}$. We define the curve $\psi:[0,1]\rightarrow\mathcal{S}_{\mathrm{nd}}$ by
\[
\psi(t) = \Psi(\Xi,\mathsf{Z}+t(\mathsf{Z}'-\mathsf{Z})),
\qquad \forall t\in[0,1].
\]
Note that $\psi$ is well-defined since $\mathcal{F}$ is convex and $\Psi$ is a diffeomorphism from $\mathrm{GL}_n\times\mathcal{F}$ to $\mathcal{S}_{\mathrm{nd}}$. Then,
\begin{align*}
\lim_{t\downarrow 0}
\frac{\pi_\gamma(\psi(t))-\pi_\gamma(\psi(0))}{t}
=\ &
\lim_{t\downarrow 0}\frac{\gamma+t(\gamma'-\gamma)-\gamma}{t} \\
=\ &
\gamma'-\gamma<0,
\end{align*}
which completes the proof.

\subsection{Proof of \Cref{proposition:main}} \label{subsection:proposition-main}
Let $\varepsilon>0$ be sufficiently small so that $\gamma'\coloneqq J(\mK')+\varepsilon<J(\mK)$. By the strict version of the bounded real lemma (see Lemma~\ref{lemma:bounded_real}), there exists $P'\succ 0$ such that
\[
\begin{bmatrix}
A_{\mathrm{cl}}(\mK')^\tr P' + P' A_{\mathrm{cl}}(\mK') & P'B_{\mathrm{cl}}(\mK') & C_{\mathrm{cl}}(\mK')^\tr \\
B_{\mathrm{cl}}(\mK')^\tr P' & - \gamma' I & D_{\mathrm{cl}}(\mK')^\tr \\
C_{\mathrm{cl}}(\mK') & D_{\mathrm{cl}}(\mK') & -\gamma' I
\end{bmatrix}\prec 0.
\]
Since the involved inequalities are strict, we can always perturb $P'$ so that $\operatorname{det}P_{12}'\neq 0$ while $P'\succ 0$ and the above inequality are still satisfied. Consequently, $(\mK',P',\gamma')\in\mathcal{S}_{\mathrm{nd}}$. Then, by the definition of $\mathcal{C}_{\mathrm{nd}}$, there exists $P\succ 0$ such that $(\mK,P,J(\mK))\in\mathcal{S}_{\mathrm{nd}}$. We can now apply~\Cref{lemma:intermediate_before_min}, which shows that there exists a $C^\infty$ curve $\psi:[0,\delta)\rightarrow\mathcal{S}_{\mathrm{nd}}$ such that $\psi(0)=(\mK,P,J(\mK))$ and
\begin{equation}\label{eq:proof_proposition4_temp}
C = \lim_{t\downarrow 0}\frac{\pi_\gamma(\psi(t))-\pi_\gamma(\psi(0))}{t} < 0.
\end{equation}
Now let
$
\varphi(t) = \pi_{\mK}(\psi(t))$ for each $t\in[0,\delta)$, where $\pi_{\mK}:\mathcal{S}_{\mathrm{nd}}\rightarrow\mathcal{C}_{\mathrm{nd}}$ denotes the canonical projection $\pi_{\mK}(\mK,P,\gamma)=\mK$. We then have $J(\mK)=\pi_\gamma(\psi(0))$ and
\begin{align*}
J(\varphi(t))
\!=\!
\inf\{\gamma\!: \exists P\!\succ\! 0\text{ s.t. }(\varphi(t),P,\gamma)\in\mathcal{S}_{\mathrm{nd}}\}
\!\leq\!
\pi_\gamma(\psi(t)).
\end{align*}
Therefore
\begin{align*}
\frac{J(\varphi(t))-J(\mK)}{t}
\leq\ &
\frac{\pi_\gamma(\psi(t))-\pi_\gamma(\psi(0))}{t}.
\end{align*}
By taking the limit superior as $t\downarrow 0$ and using~\eqref{eq:proof_proposition4_temp}, we get
\begin{equation}\label{eq:limsup_J_strictly_negative_proof}
\begin{aligned}
\limsup_{t\downarrow 0}
\frac{J(\varphi(t)) \!-\! J(\mK)}{t}
\leq\ &
C<0.
\end{aligned}
\end{equation}

Now let $\mathsf{V}=\varphi'(0)$.
we shall show that $\mathsf{V}\neq 0$, and that the left-hand side in~\eqref{eq:limsup_J_strictly_negative_proof} is in fact equal to the ordinary directional derivative in the direction $\mathsf{V}$. By~\Cref{lemma:H_inf_some_local_Lipschitz}, there exist $\epsilon>0$ and $L>0$ such that
$
|J(\mK_1)-J(\mK_2)|\leq L \|\mK_1-\mK_2\|
$
for any $\mK_i$ with $\|\mK_i-\mK\|\leq\epsilon$ for $i=1,2$. Since $\varphi$ is a $C^\infty$ curve with $\varphi(0)=\mK$, we can, without loss of generality, pick $\delta$ to be sufficiently small so that $\|\varphi(t)-\mK\|<\epsilon$ for all $t\in[0,\delta)$. Furthermore, since $\phi(t)$ is a $C^\infty$ curve, we can find $M>0$ such that
\[
\|\varphi(t)-\mK-t\,\varphi'(0)\| \leq \frac{1}{2}Mt^2,
\qquad\forall t\in[0,\delta).
\]
If $\mathsf{V}=\varphi'(0)=0$, we then have
$
|J(\varphi(t))-J(\mK)|\leq 
L\|\varphi(t)-\mK\|\leq LM t^2/2$,
which would imply
\[
\lim_{t\downarrow 0}\left|\frac{J(\varphi(t))-J(\mK)}{t}\right|=0,
\]
contradicting~\eqref{eq:limsup_J_strictly_negative_proof}. Therefore we can conclude that $\mathsf{V}\neq 0$. Finally, notice that for $t\in[0,\delta)$, we have
\begin{align*}
& \frac{J(\mK+t\,\varphi'(0))-J(\mK)}{t} \\
\leq\ &
\frac{J(\varphi(t))-J(\mK)}{t}
+\frac{|J(\mK+t\,\varphi'(0))-J(\varphi(t))|}{t} \\
\leq\ &
\frac{J(\varphi(t))-J(\mK)}{t}
+\frac{L\|\mK+t\,\varphi'(0)-\varphi(t)\|}{t} \\
\leq\ &
\frac{J(\phi(t))-J(\mK)}{t}
+\frac{1}{2}Mt.
\end{align*}
By taking the limit superior as $t\downarrow 0$ and noting that the directional derivative of $J$ always exists, we see that
\begin{align*}
J'(\mK,\mathsf{V})
=\ &
\lim_{t\downarrow 0}
\frac{J(\mK+t\,\varphi'(0))-J(\mK)}{t} \\
\leq\ &
\limsup_{t\downarrow 0}
\frac{J(\varphi(t))-J(\mK)}{t}
<0,
\end{align*}
and we arrive at the desired conclusion.

\section{Numerical Experiment}\label{sec:examples}

In this section, we provide some numerical evidence suggesting that the set $\mathcal{C}\backslash\mathcal{C}_{\mathrm{nd}}$ has measure zero.

We consider the $\mathcal{H}_\infty$ control problem for the LTI system
\begin{equation}\label{eq:LTI_numerical}
\begin{aligned}
\dot{x}(t) =\ & -x(t) + \begin{bmatrix}
    1 & 0
\end{bmatrix} w(t) + u(t), \\
z(t) =\ & \begin{bmatrix} x(t) \\ u(t) \end{bmatrix}, 
\qquad y(t) = x(t) + \begin{bmatrix} 0 & 1\end{bmatrix} w(t),
\end{aligned}
\end{equation}
where $x(t),u(t),y(t)\in\mathbb{R}$ and $z(t),w(t)\in\mathbb{R}^2$. The dynamic controller will then be parameterized by $\mK=\begin{bmatrix}
    D_\mK & C_\mK \\ B_\mK & A_\mK
\end{bmatrix}\in\mathbb{R}^{2\times 2}$. Our task is to numerically search for points in $\mathcal{C}\backslash\mathcal{C}_{\mathrm{nd}}$, and inspect whether they form a set of measure zero. Note that dynamic controllers with the same value of $B_\mK C_\mK$ will be similarity transformations of each other. Therefore, for visualization purposes, we fix $C_\mK=1$ and only examine the set $\{\mK\in\mathcal{C}:C_\mK=1\}$ instead. We also impose the constraints $
A_\mK\in[-2,2], B_{\mK}\in[-4,4],D_{\mK}\in[-1.5,1.5]$ when searching over the set $\{\mK\in\mathcal{C}:C_\mK=1\}$.

We first generate a set of points $\{\mK_j\}_{j=1}^N$ by discretizing the region $[-2,2]\times[-4,4]\times[-1.5,1.5]$ into a spatial grid with $N=101\times101\times 61$ points that are equally spaced. Then for each $j=1,\ldots,N$, we numerically compute $\gamma_j=J(\mK_j)$, and try to construct $P_j\succeq 0$ such that $\mathscr{N}(\mK_j,P_j,\gamma_j)\preceq 0$.\footnote{
Due to numerical errors, we can only find an approximate value $\hat{\gamma}_j$ of $J(\mK_j)$. In our numerical experiments, we set the tolerance so that $|\hat{\gamma}_j-J(\mK_j)|/J(\mK_j)<\epsilon$ and find $P_j$ satisfying $\mathscr{N}(\mK_j,P_j,\hat{\gamma}_j/(1-\epsilon))\preceq 0$ instead, where $\epsilon=10^{-9}$. We employ the Riccati-equation-based approach for finding $P_j$ when the associated Riccati equation is well-posed and has a positive definite solution, and turn to the LMI-based approach if the Riccati-equation-based approach does not work.
}
We then check whether the minimum eigenvalue of $P_j$ is sufficiently bounded away from zero (say greater than or equal to $10^{-4}$), and record the value of $(P_j)_{12}$.

\begin{figure}
\vspace{3pt}
    \centering
    \includegraphics[width=.45\textwidth]{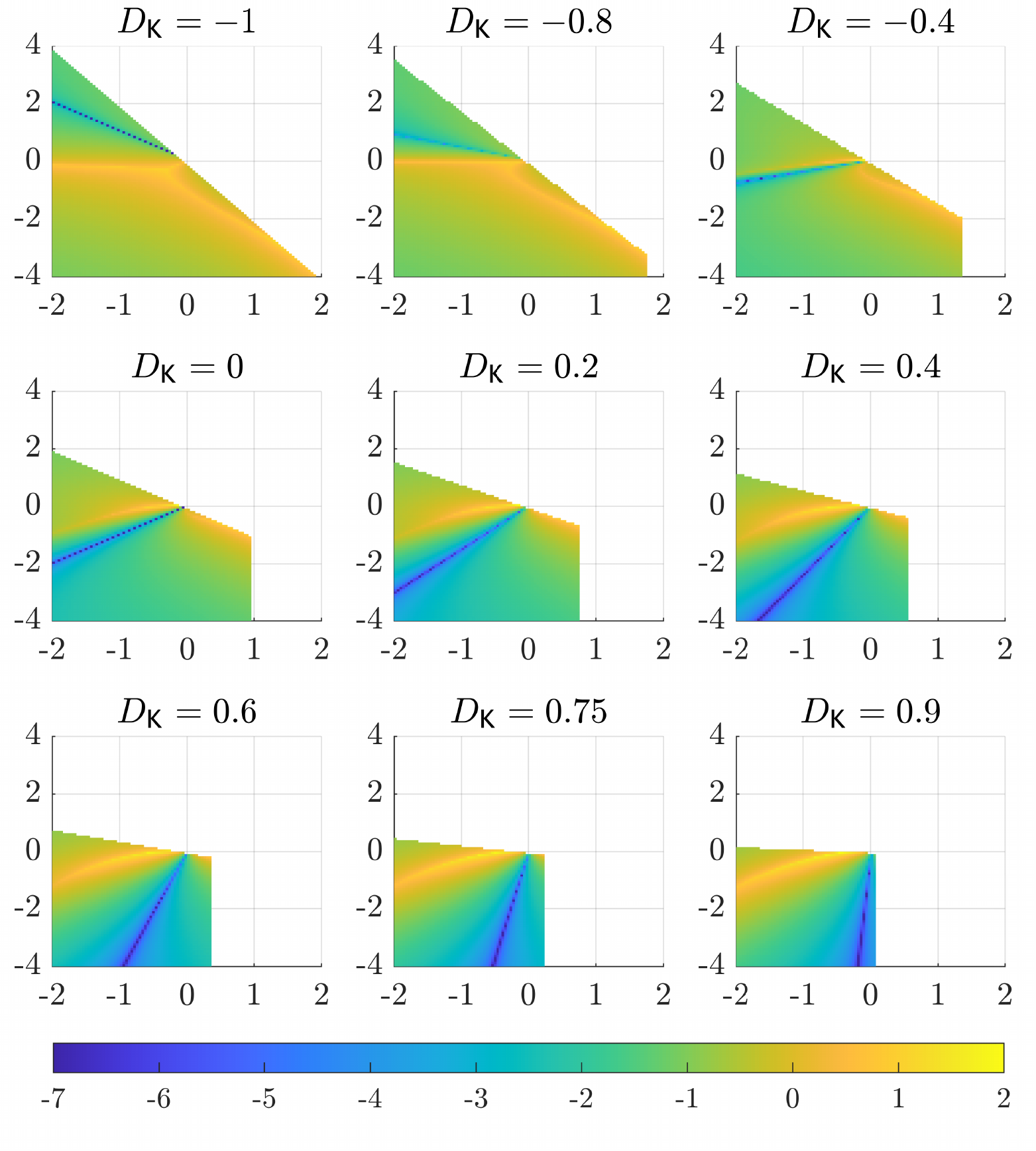}
    \caption{\footnotesize Heatmaps of $\ln|P_{12}|$ for different values of $D_{\mathsf{K}}$. The $x$-axes and $y$-axes represent $A_{\mK}$ and $B_{\mK}$ respectively. Points with very low values of $\ln|P_{12}|$ (i.e., points whose corresponding $P_{12}$ are very close to $0$) are colored in dark blue, and we can observe that they roughly form a line passing through $(0,0)$ in each sub-figure.}
    \label{fig:sim_results}
    \vspace{-12pt}
\end{figure}

Our numerical experiments show that we can find matrices $P_j$ satisfying $P_j\succ 0$ and $\mathscr{N}(\mK_j,P_j,\gamma_j)\preceq 0$ for all $j$ in the test case. \Cref{fig:sim_results} illustrates several typical heatmaps of $\ln|P_{12}|$ with fixed $D_\mK$ and varying $(A_{\mK},B_{\mK})$, generated from the recorded values $\{\ln|(P_{j})_{12}|\}_{j=1}^N$. It can be observed from the heatmaps that for each fixed value of $D_{\mK}$, the points with very low values of $\ln|P_{12}|$ seem to lie near a straight line that passes through $(0,0)$. These observations seem to suggest that, for the LTI system~\eqref{eq:LTI_numerical}, the points in $\mathcal{C}\backslash\mathcal{C}_{\mathrm{nd}}$ with $C_\mK=1$ and some fixed $D_\mK$ form a straight line passing through $(0,0)$ with a slope depending on $D_\mK$, and consequently, the set $\mathcal{C}\backslash\mathcal{C}_{\mathrm{nd}}$ could be represented as
\begin{align*}
\!\!\left\{
\!\begin{bmatrix}
D_\mK & \!\!\!\!C_\mK \\ B_\mK & \!\!\!\!A_\mK
\end{bmatrix}\!\!\in\!\mathcal{C}\! :
\cos\theta(D_\mK) \!\cdot\! A_\mK
+\sin\theta(D_\mK) \!\cdot\! B_\mK C_\mK = 0
\!\right\}
\end{align*}
for some function $\theta(D_\mK)$ of $D_\mK$, which has measure zero.

We remark that the above claim is only based on numerical results but not on rigorous derivation. Nevertheless, we believe that such results can indeed serve as numerical evidence supporting the conjecture that $\mathcal{C}\backslash\mathcal{C}_{\mathrm{nd}}$ has measure zero. The code can be found at {\small \url{https://github.com/tyj518/H_inf_Global_Optimality}}.

\section{Conclusions} \label{section:conclusion}

We consider the policy optimization for output-feedback $\mathcal{H}_\infty$ control and show that the class of \emph{non-degenerate} Clarke stationary points are globally optimal controllers, providing a global optimality certificate for direct policy search methods.
Future directions include examining whether $\mathcal{C}\backslash\mathcal{C}_{\mathrm{nd}}$ has measure zero, designing data-driven approaches for checking whether a controller is non-degenerate, convergence analysis of model-free policy search methods for $\mathcal{H}_\infty$ control, etc.

%%%%%%%%%%%%%%%%%%%%%%%%%%%%%%%%%%%%%%%%%%%%%%%%%%%%%%%%%%%%%%%%%%%%%%%%%%%%%%%%
\section*{Appendix}

\subsection{Proof of \Cref{lemma:H_inf_some_local_Lipschitz}} \label{appendix:Hinf-cost-function}

Our proof will follow the idea sketched in~\cite[Proposition 3.1]{apkarian2006nonsmooth2}, i.e., the subdifferential regularity of $\|\mathbf{T}_{zw}\|_\infty$ follows from the convexity of $\|\cdot\|_\infty$ and the continuous differentiability of the mapping from $\mathcal{C}$ to $\mathcal{RH}_\infty$ given by $\mK\mapsto \mathbf{T}_{zw}$ in \eqref{eq:transfer-function-Tzw}. But we will fill in the missing details of why the mapping $\mK\mapsto \mathbf{T}_{zw}$ is continuously differentiable. We will temporarily use $n$ to denote the dimension of $A_{\mathrm{cl}}(\mK)$, and denote $\mathcal{A} = \{A\in\mathbb{R}^{n\times n}:A\text{ is stable}\}$.

We first define the mapping $\mathscr{T}:\mathcal{A}\rightarrow\mathcal{RH}_\infty$ by
\[
(\mathscr{T}(A))(s) = (sI-A)^{-1}.
\]
Note that $A_{\mathrm{cl}}(\mK)$, $B_{\mathrm{cl}}(\mK)$, $C_{\mathrm{cl}}(\mK)$ and $D_{\mathrm{cl}}(\mK)$ are all affine functions of $\mK$ and thus are continuously differentiable. As a result, the continuous differentiability of $\mK\mapsto \mathbf{T}_{zw}$ will follow if we can show that $\mathscr{T}$ is continuously differentiable.

Let $A\in\mathcal{A}$ be an arbitrary stable matrix, and define the linear mapping $\xi_A:\mathbb{R}^{n\times n}\rightarrow \mathcal{RH}_\infty$ by
\[
(\xi_A(\Delta))(s) = (sI-A)^{-1}\Delta(sI-A)^{-1}.
\]
Let $\Delta\in\mathbb{R}^{n\times n}$ be an arbitrary matrix satisfying
$
 \|(sI-A)^{-1}\|_\infty\cdot\|\Delta\|_2<1
$, and we consider bounding the quantity
\begin{align*}
r_A(\Delta)\coloneqq\ &
\big\|\mathscr{T}(A+\Delta)-\mathscr{T}(A)-\xi_A(\Delta)
\big\|_\infty.
\end{align*}
Then, as long as $\frac{r_A(\Delta)}{\|\Delta\|_2}\rightarrow 0$ as $\|\Delta\|\rightarrow 0$, we can conclude that $\xi_A$ is the Fr\'{e}chet derivative of $\mathscr{T}$ at $A$. Indeed, we have, for any $s=\mathrm{j}\omega$ with $\omega\in\mathbb{R}$,
\begin{align*}
&(\mathscr{T}(A+\Delta)-\mathscr{T}(A)-\xi_A(\Delta))(s) \\
= &
\left((I-(sI\!-\!A)^{-1}\Delta)^{-1}-I-(sI\!-\!A)^{-1}\Delta\right)(sI\!-\!A)^{-1}.
\end{align*}
Since $\|(sI\!-\!A)^{-1}\Delta\|_\infty\leq \|(sI\!-\!A)^{-1}\|_\infty\|\Delta\|_2< 1$,
\begin{align*}
(I-(sI-A)^{-1}\Delta)^{-1}
=
\sum\nolimits_{k=0}^\infty \left((sI-A)^{-1}\Delta\right)^k,
\end{align*}
where the right-hand side converges absolutely. Therefore
\begin{align*}
r_A(\Delta) =\ 
&
\left\|\sum\nolimits_{k=2}^\infty \left((sI-A)^{-1}\Delta\right)^k (sI-A)^{-1}
\right\|_\infty \\
\leq\ &
\sum\nolimits_{k=2}^\infty (\left\|(sI \!-\!A)^{-1}\right\|_\infty \|\Delta\|_2)^k \|(sI \!-\!A)^{-1}\|_\infty \\
=\ &
\frac{\|(sI-A)^{-1}\|_\infty^3 }{1-\|(sI-A)^{-1}\|_\infty \|\Delta\|_2} \|\Delta\|^2_2,
\end{align*}
from which we can easily check that $r_A(\Delta)/\|\Delta\|_2$ converges to $0$ as $\|\Delta\|_2\rightarrow 0
$. Therefore $\mathscr{T}$ is differentiable at $A$, and its Fr\'{e}chet derivative is given by the linear mapping $\xi_A$.

Finally, we show that $\xi_A$ is continuous in $A$. By definition, we have $\xi_A(\Delta) = 
\mathscr{T}(A)\Delta\mathscr{T}(A)$, which leads to
\begin{align*}
& \sup_{\|\Delta\|_2=1}
\left\|\xi_{A+H}(\Delta)-\xi_{A}(\Delta)\right\|_\infty \\
=\ &
\sup_{\|\Delta\|_2=1}
\left\|\mathscr{T}(A+H)\Delta\mathscr{T}(A+H)
-\mathscr{T}(A)\Delta\mathscr{T}(A)\right\|_\infty \\
\leq\ &
\sup_{\|\Delta\|_2=1}
\big(\| (\mathscr{T}(A+H)-\mathscr{T}(A))\Delta \mathscr{T}(A+H)\|_\infty \\
&\qquad\qquad + \| \mathscr{T}(A)\Delta (\mathscr{T}(A+H)-\mathscr{T}(A))\|_\infty\big) \\
\leq\ &
\| (\mathscr{T}(A\!+\!H) \!-\! \mathscr{T}(A))\|_\infty
\left(\|\mathscr{T}(A\!+\!H)\|_\infty \!+\! \|\mathscr{T}(A)\|_\infty\right).
\end{align*}
As $\|H\|\rightarrow 0$, the quantity on the left-hand side will then converge to $0$, implying that the mapping $\xi_A$ is continuous in $A$. Our proof is now complete.

%%%%%%%%%%%%%%%%%%%%%%%%%%%%%%%%%%%%%%%%%%%%%%%%%%%%%%%%%%%%%%%%%%%%%%%%%%%%%%%%

\bibliographystyle{IEEEtran}
\bibliography{references}

\end{document}